\documentclass{siamltex}

\usepackage{ucs}
\usepackage[utf8x]{inputenc}
\usepackage{amsmath}
\usepackage{cite}

\usepackage{fontenc}
\usepackage[pdftex]{graphicx}

\usepackage[pdftex,
            pdfauthor={Wolfgang Bangerth, Timo Heister, Luca Heltai, Guido Kanschat, Martin Kronbichler, Matthias Maier, Toby D. Young},
            pdftitle={The deal.II Library, Version 8.1, 2013},
            colorlinks,linkcolor=blue,urlcolor=blue,citecolor=blue]{hyperref}

\newcommand{\specialword}[1]{\texttt{#1}}
\newcommand{\dealii}{{\specialword{deal.II}}}

\newcommand{\cmake}{{\specialword{CMake}}}

\title{The \dealii{} Library, Version 8.1}
\author{
  Wolfgang Bangerth%
  \thanks{Department of Mathematics, Texas A\&M University, College Station, 
    TX 77843, USA,
    \href{mailto:bangerth@math.tamu.edu}{\texttt{bangerth@math.tamu.edu}}}
  \and
  Timo Heister%
  \thanks{Mathematical Sciences,
O-110 Martin Hall.
Clemson University.
Clemson, SC 29634, USA,
    \href{mailto:heister@clemson.edu}{\texttt{heister@clemson.edu}}}
  \and
  Luca Heltai%
  \thanks{SISSA - International School for Advanced Studies, Via
    Bonomea 265, 34136 Trieste, Italy,
    \href{mailto:luca.heltai@sissa.it}{\texttt{luca.heltai@sissa.it}}}
  \and
  Guido Kanschat%
  \thanks{Interdisciplinary Center for Scientific Computing, Heidelberg
          University, Im Neuenheimer Feld 368, 69120 Heidelberg, Germany,
          \href{mailto:kanschat@uni-heidelberg.de}{\texttt{kanschat@uni-heidelberg.de}}}
  \and
  Martin Kronbichler%
  \thanks{Institute for Computational Mechanics, Technische Universit\"at M\"unchen, Boltzmannstr.~15, 85748 Garching b. M\"unchen, Germany, \href{mailto:kronbichler@lnm.mw.tum.de}{\texttt{kronbichler@lnm.mw.tum.de}}}
  \and
  Matthias Maier%
  \thanks{Institute of Applied Mathematics, Heidelberg
          University, Im Neuenheimer Feld 293/294, 69120 Heidelberg, Germany,
          \href{mailto:matthias.maier@iwr.uni-heidelberg.de}{\texttt{matthias.maier@iwr.uni-heidelberg.de}}}
  \and
  Bruno Turcksin%
  \thanks{Department of Mathematics, Texas A\&M University, College Station, 
    TX 77843, USA,
    \href{mailto:turcksin@math.tamu.edu}{\texttt{turcksin@math.tamu.edu}}}
  \and
  Toby~D.~Young%
  \thanks{Institute of Fundamental Technological Research of the Polish Academy of Sciences, ul. Pawi\'nskiego 5b, Warsaw 02-106, Poland,
    \href{mailto:tyoung@ippt.pan.pl}{\texttt{tyoung@ippt.pan.pl}}}
}

\begin{document}
\maketitle

\begin{abstract}
  This paper provides an overview of the new features of the finite element library \dealii{} version 8.1.
\end{abstract}

\section{Overview}

\dealii{} version 8.1 was released December 24, 2013. This paper provides an
overview of the new features of this release and serves as a citable
reference for the \dealii{} software library version 8.1. \dealii{} is an
object-oriented finite element library used around the world in the
development of finite element solvers. It is available for free under the
GNU Lesser General Public License (LGPL) from the \dealii{} homepage at
\url{http://www.dealii.org/}.

Version 8.1 contains,  along
with the usual set of bug fixes and documentation updates, the following noteworthy changes:
\begin{itemize}
\item Three new tutorial programs (see Section \ref{sec:step26}, \ref{sec:step42}, and \ref{sec:step51});
\item Improved support for multicore parallelization on shared memory machines (Section \ref{sec:threading});
\item The testsuite was ported to CTest/CDash (Section \ref{sec:testsuite});
\item Post-install tests (see Section \ref{sec:installtests}).
\end{itemize}
Information on how to cite \dealii{} is provided in Section \ref{sec:cite}.

\section{Changes to the governance structure}

\dealii{} has been a project with contributors from around the world for a
long time already. However, there has never been a formalized way to recognize
contributions other than by listing the authors in a file that is part of
the documentation.

Starting with this release, \dealii{} now has a more open governance structure
that we hope will more accurately reflect the extensive contributions of many
participants in this project: in addition to listing all contributors as
before, \dealii{} is now governed by a council of \textit{developers} --
currently Luca Heltai, Martin Kronbichler, Matthias Maier, Bruno Turcksin, and
Toby Young -- along with \textit{principal developers} who are responsible for
running the technical infrastructure -- currently Wolfgang
Bangerth, Timo Heister, and Guido Kanschat. It is our hope that this new
structure is a more adequate representation of many contributors' long-term
effort for this project. It is our intention that this structure remains open
to others and that this serves as motivation for others to participate!

\section{Changes to the library}

\subsection{A new tutorial: Elasto-plastic contact problem (step-42)}
\label{sec:step42}

The first of the three new tutorial programs, step-42 (written by Jörg Frohne,
Timo Heister, and Wolfgang Bangerth), shows how to solve an elasto-plastic
contact problem. The program is an extension of step-41 to a much more complex
equation (nonlinear elasto-plasticity) and also demonstrates how to compute
solutions for such problems in parallel.

The example shows how to solve a 3d, non-linear problem with a semi-smooth
Newton method coupled with an active set strategy on adaptively refined meshes
and scales well up to at least 1024 cores and 1 billion unknowns. An
accompanying paper \cite{FHB13} explains the details and design of the
algorithms behind this tutorial. It also shows scalability results for
parallel computations.

\subsection{A new tutorial: A hybridizable discontinuous Galerkin method (step-51)}
\label{sec:step51}

Step-51, written by Martin Kronbichler and Scott Miller, shows how to
implement a hybridizable discontinuous Galerkin method (HDG) in \dealii{}. An
HDG method is a special DG method which features a reduced number of
globally coupled degrees of freedom compared to other DG schemes for implicit
systems. This is achieved by
introducing a new variable for the numerical trace, i.e., the operator that
connects the solution of the subproblems on each element. The linear system to
be solved globally then only consists of degrees of freedom in the trace variable,
whereas the DG contributions interior to cells are eliminated during assembly by static
condensation. The tutorial program shows how to implement this concept in \dealii{},
where the class \texttt{FE\_FaceQ} represents the trace element for the usual
DG element \texttt{FE\_DGQ} applied to the scalar convection-diffusion
equation written as a first order system. 

The tutorial program includes a few
practical aspects of an HDG implementation, namely the realization of
superconvergent post-processing and the utilization of the parallel character
of the assembly and static condensation process with the \texttt{WorkStream}
class. The tutorial also contains an extensive discussion of efficiency of
this approach compared to standard continuous finite elements. The results
show that the HDG method is competitive to continuous elements for medium
order between about 3 and 6 in terms of solution time on diffusion-dominated
systems. For convection-dominated problems, the method inherits the superior
stability properties of low to medium order DG methods, which allows for good
solution quality already on coarser meshes.

Furthermore, a new element \texttt{FE\_FaceP} has been added to the \dealii{}
library that can used for hybridization of the discontinuous elements of
complete degree $p$, \texttt{FE\_DGP}. \texttt{FE\_DGP} and \texttt{FE\_FaceP}
of degree $p$ have fewer degrees of freedom per cell compared to the elements
\texttt{FE\_DGQ} and \texttt{FE\_FaceQ} of tensor product degree $p$.

\subsection{A new tutorial: An adaptive mesh solver for the heat
  equation (step-26)}
\label{sec:step26}

\dealii{} has long had tutorial programs for the wave equation, but none for
the conceptually simpler heat equation. Furthermore, there was no simple
tutorial program that showed how one can solve time dependent problems where
the mesh changes every few time steps.

The new step-26 tutorial program closes this gap. It is a rather simple
program -- the program has just 151 lines with a semicolon -- and as such
serves as a gentle introduction to these topics.

\subsection{Improved WorkStream implementation}
\label{sec:threading}
The \texttt{WorkStream} namespace has already in the past contained the
functions that are used to parallelize (using multithreading) many of the loops over all cells or
faces that one encounters in finite element computations. Extensive
documentation on how to use this class and the rationale behind it is provided
in the ``Parallel computing with multiple processors accessing shared memory''
module.

This class has been significantly revamped in an effort to improve scaling to
large multicore workstations. In particular, we now use thread-local variables
in some places to improve cache performance. In addition, the namespace has
obtained another implementation of \texttt{WorkStream::run} that does not just
take a sequence of cells (or other objects to work on), but such a sequence
that has been colored in a way to indicate which cells will conflict with
others when writing into the global matrix, right hand side, or other
object. Using this function instead of the previously existing one obviously
requires changes to existing code bases, but can provide significant speedups
in some circumstances when using significant numbers of threads (say, above 16
or 32). The details of these algorithms have been documented in \cite{TKB13}.

As part of the changes necessary to allow this, \dealii{} has also gained a
set of generic functions that provide graph coloring algorithms. The are
located in namespace \texttt{GraphColoring}.

\subsection{Testsuite ported to CTest/CDash}
\label{sec:testsuite}

With this minor release \dealii{} is now fully ported to \cmake. The last
remaining step was the migration of the testsuite to \textsc{CTest} as test
driver and \textsc{CDash} as web front end. This involved porting over 3000
tests to a new directory structure. Furthermore, the regression and build
tests are no longer independent testsuites, but a combined one; and tests
are by default run against both the debug and release versions of the
library.

The port of the testsuite was also motivated by the fact that, in order to
provide official support for at least three major compiler brands
(\textsc{GCC}, \textsc{Clang} and \textsc{ICC}) on multiple platforms such
as Linux, MAC OSX and BSD, it is highly necessary to have regular
regression tests available for these platforms. Hence, in addition to the
regression tester already present (whose main purpose is to test every
single subversion revision) dedicated nightly build tests for a big variety
of above compilers and platforms (and external dependencies) have been set
up.

The testsuite is set up in the build directory via \texttt{make
setup\_tests}. After that, tests can be run by invoking the test driver
\texttt{ctest} with suitable options. In order to submit test results a
\texttt{CTest} script is available that can be passed as an option to the
driver, e.\,g. by \texttt{ctest -S ../tests/run\_testsuite.cmake} (or if
just the build tests should be run, \texttt{\dots/run\_buildtest.cmake}).
The scripts run a configure, build and test stage as necessary and submit
the results to a \texttt{CDash} instance. Public, anonymous submission
of test results to the \textsc{CDash} site is possible in order to
easily upload and share regression test results.

\subsection{Post-install tests}
\label{sec:installtests}

Every installation now ships with a small collection of tests that can be executed by calling \texttt{make test} in the build directory. This is used to identify common configuration problems, bugs in external dependencies, problems in the MPI wrappers, and to check the correct setup of some of the external packages.
Running the tests is now part of the installation instructions and we especially recommend to run the tests on new machines.

\subsection{Other Changes}\label{sec:other}

The \dealii{} release 8.1 also includes improvements in the following areas:

\begin{itemize}
\item The handling of periodic boundary conditions has been improved and extended
  to distributed triangulations.
\item \dealii{}'s output capabilities have been extended in several
  aspects. In particular, it is now possible to merge output from vectors
  belonging to several different \texttt{DoFHandler} objects (a common case
  when solving multiphysics problems) within a single file instead of forcing
  the user to join the data in one \texttt{DoFHandler} for output.
\item Higher order polynomial boundary descriptions now use support points of
  the Gauss-Lobatto quadrature formula instead of equidistant ones, which
  gives much better stability at high orders and more accuracy for some curved
  boundary shapes due to superconvergence effects.
\end{itemize}

\section{How to cite \dealii{}}\label{sec:cite}

In order to justify the work the developers of \dealii{} put into this
software, we ask that papers using the library reference one of the \dealii{}
papers. This helps us justify the effort we put into it.

There are various ways to reference \dealii{}. To acknowledge the use of a
particular version of the library, reference the present document as follows:
\begin{verbatim}
@article{dealII81,
  title = {The {\tt deal.{I}{I}} Library, Version 8.1},
  author = {W. Bangerth and T. Heister and L. Heltai and G. Kanschat
   and M. Kronbichler and M. Maier and B. Turcksin and T. D. Young},
  journal = {arXiv preprint \url{http://arxiv.org/abs/1312.2266v4}},
  year = {2013},
}
\end{verbatim}

The original \texttt{\dealii{}} paper containing an overview of its
architecture is \cite{BangerthHartmannKanschat2007}. If you rely on specific
features of the library, please consider citing any of the following:
\begin{itemize}
 \item For geometric multigrid: \cite{Kanschat2004,JanssenKanschat2011};
 \item For distributed parallel computing: \cite{BangerthBursteddeHeisterKronbichler11};
 \item For $hp$ adaptivity: \cite{BangerthKayserHerold2007};
 \item For matrix-free and fast assembly techniques:
   \cite{KronbichlerKormann2012};
 \item For computations on lower-dimensional manifolds:
   \cite{DeSimoneHeltaiManigrasso2009}.
\end{itemize}

\nocite{BangerthKanschat1999}

\section{Acknowledgements}

\dealii{} is a world-wide project with dozens of contributors around the
globe. Other than the authors of this paper, the following people contributed code to
this release: 
  Fahad Alrashed,
  Daniel Arndt,
  Juan Carlos Araujo Cabarcas,
  Krzyszof Bzowski,
  Francesco Cattoglio,
  Denis Davydov,
  David Emerson,
  Armin Ghajar Jazi,
  Eric Heien,
  Tobin Isaac, 
  Oleh Krehel,
  Craig Michoski,
  Scott Miller,
  Jean-Paul Pelteret,
  Andreas Putz,
  Mayank Sabharwal,
  Martin Steigemann.
Their contributions are much appreciated!

\dealii{} and its developers are financially supported through a
variety of funding sources. W.~Bangerth and B.~Turcksin were partially
supported by the National Science Foundation under award OCI-1148116
as part of the Software Infrastructure for Sustained Innovation (SI2)
program; by the Computational Infrastructure in Geodynamics initiative
(CIG), through the National Science Foundation under Award
No.~EAR-0949446 and The University of California -- Davis; and through
Award No.~KUS-C1-016-04, made by King Abdullah University of Science
and Technology (KAUST). 

L.~Heltai was partially supported by the project OpenSHIP,
``Simulazioni di fluidodinamica computazionale (CFD) di alta qualit\`a
per le previsioni di prestazioni idrodinamiche del sistema
carena-elica in ambiente OpenSOURCE'', financed by Regione FVG - POR
FESR 2007–2013 Obiettivo competitivit\`a regionale e occupazione.

T.~Heister was partially supported by the Computational Infrastructure in
Geodynamics initiative (CIG), through the National Science Foundation under Award
No. EAR-0949446 and The University of California -- Davis and through Award No. KUS-
C1-016-04, made by King Abdullah University of Science and Technology (KAUST).

The Interdisciplinary Center for Scientific Computing (IWR) at Heidelberg University has provided 
hosting services for the deal.II web page and the SVN archive.
 

\bibliography{deal81}{}
\bibliographystyle{plain}

\end{document}